\documentclass[11pt]{amsart}

\usepackage{verbatim, amscd, epsfig,amsmath,amsthm,amstext,amssymb,hyperref
}
\theoremstyle{plain}

\input xypic

\newtheorem{theorem}{Theorem}[section]

\theoremstyle{definition}

\theoremstyle{remark}
\newtheorem{rem}[theorem]{Remark}

\newcommand{\del}{\partial}
\newcommand{\delbar}{\bar{\del}}

\newcommand{\R}{\mathbb{ R}}
\newcommand{\C}{\mathbb{ C}}

\renewcommand{\H}{\mathbb{ H}}

\renewcommand{\P}{\mathbb{ P}}
\newcommand{\HP}{\H\P}
\newcommand{\CP}{\C\P}

\DeclareMathOperator{\End}{End}
\DeclareMathOperator{\Hom}{Hom}

\DeclareMathOperator{\im}{im}

\begin{document}

\title
{Sequences of Willmore surfaces}
\author{K. Leschke and F. Pedit}
\address{
Katrin Leschke\\
Institut f\"ur Mathematik\\
Lehrstuhl f\"ur Analysis und Geometrie\\
Universit\"at Augsburg\\
86135 Augsburg\\
Germany
}
\address{Franz Pedit\\
Department of Mathematics and Statistics \\
University of Massachusetts\\
Amherst, MA 01003, USA\\
and\\
Mathematisches Institut\\
Universit\"at T\"ubingen\\
Auf der Morgenstelle 10\\ 
72076 T\"ubingen\\
Germany}
\email{katrin.leschke@math.uni-augsburg.de, pedit@math.umass.edu}

\begin{abstract}In this paper we develop the theory of Willmore
  sequences for Willmore surfaces in the $4$-sphere. We show that
  under appropriate conditions this sequence has to terminate. In this
  case the Willmore surface either is the twistor projection of a
  holomorphic curve into $\CP^3$ or the inversion of a minimal surface
  with planar ends in $\R^4$. These results give a unified explanation
  of previous work on the characterization of Willmore spheres and
  Willmore tori with non-trivial normal bundles by various authors.
\end{abstract}

\thanks{The first author thanks the Department of Mathematics and Statistics at the University of Massachusetts, Amherst, and the Center for Geometry, Analysis, Numerics and Graphics for their support and hospitality. The second author thanks the Humboldt Foundation and the Matheon Center at TU-Berlin for their support}
\maketitle
\section{introduction}
The differential geometric transformation theory of special surfaces
classes plays an important role in the study of their classification:
initially, these transformations were used to construct more
complicated examples of surfaces from ``trivial'' surfaces and, more
recently, the existence of such transformations has been linked to the
phenomenon of complete integrability. Classically known examples
include the B\"acklund transformations of constant Gaussian or mean
curvature surfaces and the Darboux transformations of isothermic
surfaces. Both types of transformations satisfy what is known as
Bianchi permutability and therefore generate an abelian group acting
on those surfaces. Under appropriate circumstances orbits of this
action can be seen as the ``energy shells'' of a completely integrable
system. In analytical terms these transformations allow the
construction of solutions of the nonlinear system of PDEs describing a
special surface class from known, often trivial, solutions by solving
auxiliary ODEs together with algebraic manipulations.

In general, of course, one cannot expect to obtain a complete
classification of all surfaces of a particular type from successive
applications of such transforms.  But there are a number of
interesting instances, most strikingly perhaps in the theory of
harmonic maps, where this is the case.  Physicists
\cite{glaser-stora},\cite{din-zak} working on non-linear sigma models
observed that taking the $(1,0)$-part of the derivative of a harmonic
map from a Riemann surface $M$ into complex projective space $\C\P^n$
yields a new harmonic map into $\C\P^n$, thereby generating what is
now called the harmonic sequence.  This construction turned out to be
important for classifying harmonic maps: tracing the energies of the
harmonic maps in the sequence one obtains, under appropriate degree
assumptions, an increasing function in the sequence length. This
function turns out to be bounded by the energy of the initial harmonic
map.  Thus the sequence has to terminate, in which case the last
element of the sequence is a holomorphic map, and the initial harmonic
map is part of the Frenet frame of a holomorphic curve in $\C\P^n$. In
other words, under an appropriate degree assumption on a harmonic map
into $\CP^n$ such a map comes from a holomorphic curve by taking
derivatives and projections \cite{durham-boys}.  This construction can
be generalized to other target spaces such as Grassmannians
\cite{wolfson} and Lie groups giving rise to the uniton transform of
Uhlenbeck \cite{uhlenbeck} .

These advances in the theory of harmonic maps from a Riemann surface
have implications to classical surface theory: in certain cases
curvature conditions on a surface translate to harmonicity conditions
on a ``Gauss map'', the classical example being that of constant mean
curvature surfaces whose unit normal maps into the $2$-sphere are
harmonic.  Applying the harmonic sequence arguments to this special
case yields the following result by Eells and Wood \cite{eells&wood}:
if the degree of a harmonic map from a compact Riemann surface $M$
into $\CP^1$ (with the appropriate orientation) is at least the genus
$g$ of $M$ then the map has to be holomorphic. For example, a constant
mean curvature sphere has a holomorphic unit normal map and thus is a
round sphere as first observed by Hopf \cite{Hopf}.  On the other
hand, immersed constant mean curvature tori have unit normal maps of
degree zero and the harmonic sequence does not terminate. This case
leads to the theory of spectral curves and the construction of
constant mean curvature tori from algebraically completely integrable
systems \cite{PS}, \cite{hitchin-harmonic}, \cite{Bob}.

In the present paper we consider Willmore surfaces in the $4$-sphere
$S^4$. These surfaces are also characterized by a harmonicity
condition, this time on their conformal Gauss maps or mean curvature
sphere congruences. Bryant's \cite{Bryant}, Ejiri's \cite{Ejiri}, and
later Montiel's \cite{montiel}, classification of Willmore spheres in
$S^3$ and $S^4$ indicate that there aught to be a harmonic sequence
type explanation for these results. But the conformal Gauss map takes
values in the space of oriented $2$-spheres in $S^4$, a Grassmannian
of space-like planes, for which the harmonic sequence and uniton
transform constructions are much less developed. To remedy this
situation and to stay conceptually close to the $\CP^n$ case, we
describe the $4$-sphere by the quaternionic projective line $\HP^1$ so
that a Willmore surface $f\colon M\to S^4$ is given by the line
sub-bundle $L\subset V$, $L_p=f(p)$, where $V=M\times\H^2$ is the
trivial bundle over $M$. The conformal Gauss map then is a harmonic
complex structure $S\in\Gamma(\End(V))$, $S^2=-1$, on $V$ with the
property that the $2$-spheres $S(p)$, given by the quaternionic
eigenlines of $S(p)$, touch the surface $f$ and have the same mean
curvature then $f$ at $p\in M$.  In terms of the type decomposition
\[
\tfrac{1}{2}SdS=A + Q
\]
of the derivative of $S$ these conditions imply 
\[
L\subset \ker Q\quad\text{or equivalently,}\quad \im A\subset L
\]
and, denoting by $*$ the complex structure on $M$,
\[
d*A=d*Q=0\,,
\]
which is the harmonicity of $S$. If $A$ and $Q$ are not identically
zero they both have rank one and we can define \cite{coimbra} two more
line sub-bundles $\hat{L}\subset V$ and $\tilde{L}\subset V$
satisfying
\[
\im Q\subset \hat{L}\,,\quad \text{and} \quad \tilde{L}\subset \ker A\,.
\]
Due to the harmonicity equations these line sub-bundles extend
smoothly across the isolated zeros of $A$ and $Q$. Moreover, since $A$
and $Q$ have type the maps $\hat{f},\tilde{f}\colon M\to S^4$ given by
the line bundles $\hat{L}$ and $\tilde{L}$ are conformal and thus have
at most isolated branch points.  Since $A$ and $Q$ anti-commute with
$S$ the bundles $\hat{L}$ and $\tilde{L}$ are stable under $S$ which
means that the surfaces $\hat{f}$ and $\tilde{f}$ lie for each $p\in
M$ on the mean curvature sphere $S(p)$ of $f$. Unless
$\hat{f}=\tilde{f}$, in which case the mean curvature sphere
congruence $S$ of $f$ has a second envelope, the mean curvature
spheres of $\hat{f}$ and $\tilde{f}$ are distinct from $S$.  By using the 
theory of $1$-step B\"acklund transforms \cite{coimbra}, \cite{osculates}, we show that
the mean curvature sphere congruences of $\hat{f}$ and
$\tilde{f}$ extend smoothly across the possible branch points. 
Since $\hat{A}=Q$ and $\tilde{Q}=A$ the conformal Gauss
maps $\hat{S}$ and $\tilde{S}$ are again harmonic and therefore
$\hat{f}$ and $\tilde{f}$ are Willmore surfaces. Continuing this
construction we obtain, starting from $f_0=f$ and putting
$f_{-1}=\hat{f}$, $f_{1}=\tilde{f}$, a sequence
\[
 \cdots\rightarrow f_{-2}\rightarrow f_{-1}\rightarrow f_0\rightarrow f_1 \rightarrow f_2\rightarrow\cdots
 \]
 of (possibly branched) Willmore surfaces $f_i\colon M\to S^4$ with
 smooth conformal Gauss maps $S_i\in\Gamma(\End(V)) $, $S_i^2=-1$.
 Not unlike in the case of the harmonic sequence,
 tracing the Willmore energies of $f_i$ along the sequence provides a
 function in the sequence length bounded by the Willmore energy of
 $f$. In case $M$ has genus zero, or $M$ has genus one and $f$ has
 non-trivial normal bundle degree, this function in the sequence
 length is strictly increasing and therefore the sequence has to
 terminate.  
 
 If the Willmore sequence terminates, it can only happen in two ways: 
 either some $A_i=0$ in which
 case already $A=0$ or $Q=0$ and $f$ is the twistor projection of a
 holomorphic curve into $\CP^3$; or one of the maps $f_i$ is constant
 in which case already $f_{-1}=f_{i}=f_{1}$ and stereographic
 projection of $f$ from this point gives a minimal surface in $\R^4$.
 Thus, if the Willmore sequence of a Willmore surface $f\colon M\to
 S^4$ terminates, $f$ is either the twistor projection of a
 holomorphic curve into $\CP^3$ or comes from a minimal surface in
 $\R^4$ via stereographic projection.
 This provides a unifying point of view of the results in
 \cite{Bryant}, \cite{Ejiri}, \cite{montiel}, and
 \cite{willmore_tori}. As already conjectured and partially verified
 in \cite{willmore_tori}, the analogue for Willmore surfaces of the
 result by Eells and Wood \cite{eells&wood} should be the following:
 if the normal bundle degree of $f$ (with the appropriate orientation)
 is at least $4g-3$, where $g$ is the genus of $M$, then the Willmore
 sequence of $f$ terminates.

\section{The Willmore sequence}
  
In the construction of the Willmore sequence of a Willmore surface in
the $4$-sphere $S^4$ one generally encounters special branched
conformal immersions $f\colon M\to S^4$ from a Riemann surface $M$
whose mean curvature spheres extend smoothly across the branch points
of $f$.  At a regular point $p\in M$ of $f$ the mean curvature sphere
is the oriented $2$-sphere $S(p)\subset S^4$ touching $f$ and having
the same mean curvature as $f$ at the point $p\in M$. Even though this
description seems to depend on Euclidean quantities it is invariant
under the M\"obius group of $S^4$. The map which assigns to $p\in M$
the mean curvature sphere $S(p)$ is called the conformal Gauss map or
the mean curvature sphere congruence and plays a pivotal role in the
theory of Willmore surfaces.

In what follows, we present the conformal $4$-sphere $S^4$ by the
quaternionic projective line $\HP^1$ on which ${\bf Gl}(2,\H)$ acts by
M\"obius transformations. A round $2$-sphere $S\subset S^4$ is given
by the quaternionic eigenlines of a complex structure $S\in {\bf
  Gl}(2,\H)$, $S^2=-1$, on $\H^2$ and $-S$ is the sphere with the
opposite orientation.  The conformal map $f\colon M\to S^4$ is
described as the line sub-bundle $L\subset V$ of the trivial bundle
$V=M\times \H^2$ over $M$, namely $L_p=f(p)$ for $p\in M$.  As always
with maps into projective spaces, under the identification
$f^*TS^4=\Hom(L,V/L)$, the derivative of $f$ is expressed by the
$1$-form
\[
\delta=\pi d_{|_L}\in\Omega^1(\Hom(L,V/L))\,.
\]
Here $d$ denotes the trivial connection on $V$ and $\pi\colon V \to
V/L$ is the projection into the quotient bundle.  A congruence of
touching $2$-spheres $S$ along the conformal map $f\colon M\to S^4$
thus is a complex structure $S$ on the bundle $V$ stabilizing
$L\subset V$ and satisfying
\[
*\delta=S\delta=\delta S\,,
\]
where $*$ denotes the complex structure on $1$-forms on $M$.  The
first equality expresses the conformality of $f$ and the second
equality assures that, at regular points $p\in M$ of $f$, the tangent
spaces to $f$ coincide with the tangent spaces of the spheres $S(p)$
at $f(p)$.  Since $S$ stabilizes $L$ there are well-defined induced
complex structures, again denoted by $S$, on $L$ and $V/L$.  The
existence of a touching sphere congruence along $f$ allows us to
decompose the pulled back tangent bundle of the $4$-sphere,
\[
f^*TS^4=\Hom(L,V/L)= \Hom_{+}(L,V/L)\oplus \Hom_{-}(L,V/L)\,,
\]
into the two complex line bundles $\Hom_{\pm}(L,V/L)$ which agree with
the tangent and normal bundles of $f$ along the regular points of
$f$. The subscript $\pm$ indicates complex linear respectively
anti-linear homomorphisms with respect to the complex structures
induced by $S$. Note that the complex structure on a $\Hom$-bundle is
always given by the target complex structure so that, for instance,
$\Hom_{-}(L,V/L)=\Hom_{+}(\bar{L},V/L)$, where $\bar{L}$ indicates
that we use the opposite complex structure $-S$ on $L$.  A
quaternionic vector bundle with a complex structure is just the double
of a complex vector bundle \cite[Sec. 11.1]{coimbra}: therefore
$V=V_{+}\oplus V_{-}$ where $V_{\pm}$ are the $\pm i$ eigen-spaces of
$S$ and $V_{-}$ is complex isomorphic to $V_{+}$ via multiplication by
$j$. Thus $\Hom_{+}(L,V/L)=\Hom_{\C}(L_+,(V/L)_{+})$ is indeed a
complex line bundle.

If $M$ is compact the degree of a quaternionic bundle with complex
structure is the degree of the underlying complex bundle: for example,
the normal bundle degree of $f$ is
\begin{equation}\label{eqn:normal-degree}
v:=\deg \Hom_{-}(L,V/L)=\deg \Hom_{+}(\bar{L},V/L)=\deg (V,S)\,.
\end{equation}
What distinguishes the mean curvature sphere congruence $S$ among
touching sphere congruences is a second order touching condition
\cite[Thm. 2]{coimbra}. In terms of the type decomposition of the
derivative $dS=dS'+dS''$ with respect to $S$ on $V$ and $*$ on $M$
this condition is
\[
L\subset \ker dS''\quad\text{or equivalently,}\quad \im dS'\subset L\,.
\]
The trivial connection $d$ on $V$ can be written in terms of the
complex connection $\hat{\nabla}$, $\hat\nabla S=0$, as
\begin{equation}\label{eqn:nabla}
d=\hat{\nabla}+\tfrac{1}{2}SdS=\hat{\nabla}+A+Q\,,
\end{equation}
where $A\in\Gamma(K\End_{-}(V))$ and $Q\in\Gamma(\bar{K}\End_{-}(V))$
are the type decomposition of $\tfrac{1}{2}SdS=A+Q$ and $K$ denotes
the canonical bundle of $M$. With this notation
\begin{equation*}\label{eqn:dS}
dS=2*(Q-A)
\end{equation*}
and the mean curvature sphere condition can be rewritten as
\[
L\subset \ker Q\quad\text{or equivalently,}\quad \im A\subset L\,.
 \]
 The Willmore energy \cite[Def. 8]{coimbra} of the conformal map
 $f\colon M\to S^4$ with mean curvature sphere congruence $S\in
 \Gamma(\End(V))$ is the same as the Dirichlet energy
\begin{equation}\label{eqn:Will-energy}
W(f)=2\int_{M} <A\wedge *A>
\end{equation}
of $S$.  At regular points of $f$ the integrant $<A\wedge *A>$ is
given by the usual Willmore integrant
$(|H|^2-K-K^{\perp})\text{vol}_{f^*h}$, where the mean curvature $H$,
Gaussian curvature $K$, and normal bundle curvature $K^{\perp}$ are
computed with respect to a conformally flat metric $h$ on $S^4$.

The conformal map $f\colon M\to S^4$ with mean curvature sphere
congruence $S$ is said to be a {\em branched Willmore surface} if $f$
is a critical point for the Willmore energy under compactly supported
variations of $f$ by maps which have a mean curvature sphere
congruence.  Note that we also allow the conformal structure of $M$ to
change. As in the unbranched case the resulting Euler-Lagrange
equation \cite{osculates} is the harmonic map equation
\begin{equation}\label{eqn:harm}
d*A=d*Q=0
\end{equation}
for $S$.  The complex linear part of the last equation gives 
\begin{equation}\label{eqn:holo}
d^{\hat{\nabla}}A=d^{\hat{\nabla}}Q=0
\end{equation}
which says that $A\in H^0(K\Hom_+(\bar{V},V))$ and $Q\in
H^0(K\Hom_+(V,\bar{V}))$ are complex holomorphic bundle maps.
Assuming that $A$ and $Q$ are not identically zero they have rank one
and we obtain two smooth line sub-bundles $\tilde{L},\hat{L}\subset
V$, namely
 \[
\tilde{L}\subset \ker A\quad\text{and}\quad \im Q\subset \hat{L}\,,
\]
extending the kernel of $A$ and the image of $Q$ across their
zeros. The corresponding smooth maps $\tilde{f},\hat{f}\colon M\to
S^4$ are called {\em B\"acklund transforms} of $f$. This construction
already appeared in \cite{coimbra} but it was unknown whether the
B\"acklund transforms admitted smooth mean curvature spheres. Their
existence is crucial for the continuation of this construction to
obtain a sequence of Willmore surfaces.
 \begin{theorem}\label{thm:main}
   Let $f\colon M\to S^4$ be a branched Willmore surface with mean
   curvature sphere congruence $S$. Then the B\"acklund transforms
   $\tilde{f},\hat{f}\colon M\to S^4$ are again branched Willmore
   surfaces with mean curvature sphere congruences $\tilde{S}$ and
   $\hat{S}$.  Moreover, $\tilde{S}=-S$ on $V/\tilde{L}$, $\hat{S}=-S$
   on $\hat{L}$, $\tilde{Q} = A$ and $\hat{A}= Q$.
 \end{theorem}
 This theorem is a consequence of a more fundamental result,
 Theorem~\ref{lem:1-stepwillmore}, about the smoothness of the mean
 curvature sphere congruences of $1$-step B\"acklund transforms. Since
 this result is interesting in its own right and relies on the theory
 of $1$-step B\"acklund transforms developed in
 \cite{coimbra},\cite{osculates},\cite{onestep}, and \cite{habil}, we
 postpone its proof to the last section of this paper.

 {F}rom Theorem~\ref{thm:main} and $\im A\subset L\subset \ker Q$ we see
 that
\[
\tilde{\hat{f}}=f=\hat{\tilde{f}}\,.
\]
Therefore we obtain a sequence, the {\em Willmore sequence},
\[
\cdots\rightarrow f_{-2}\rightarrow f_{-1}\rightarrow f_0=f\rightarrow f_{1}\rightarrow f_{2}\rightarrow\cdots
\]
of branched Willmore surfaces $f_i\colon M\to S^4$,
$f_{i+1}=\tilde{f_i}$, with mean curvature spheres $S_i$ provided that
the $A_i$ do not vanish and that the $f_i$ are not constant. The
difference of the Willmore energies of successive sequence elements is
given by the normal bundle degree \eqref{eqn:normal-degree} of $f_i$,
\begin{equation}\label{eqn:Will-diff}
W(f_i)-W(f_{i-1})=4\pi\deg(V,S_i)=4\pi v_i\,.
\end{equation}
This follows from
\[
W(f)-W(\hat{f})=2\int_M <A\wedge *A> - 2\int_M <Q\wedge *Q>=2\int_M S R^{\hat{\nabla}}=4\pi\deg(V,S)=4\pi v\,,
\]
where we used \eqref{eqn:Will-energy} and the complex linear part of
the zero curvature equation for $d= \hat{\nabla}+A+Q$.

We are now in the position to give a variation of the arguments used
in the harmonic sequence construction
\cite{wolfson},\cite{durham-boys}. Eventually, this will give a
characterization of Willmore surfaces under some assumptions on their
normal bundle degrees.  First we note that $*\delta=S\delta=\delta S$
implies that the holomorphic structure $\delbar=\hat{\nabla}''$
stabilizes $L$ and that $\delbar \delta=0$. This means that the
derivative $\delta$ of $f\colon M\to S^4$ is a holomorphic section of
the complex line bundle $K\Hom_+(L,V/L)$. We already have seen that
$A\in H^0(K\Hom_+(\bar{V},V))$ is a complex holomorphic bundle map
\eqref{eqn:holo}.  But the anti-holomorphic structure
$\del=\hat{\nabla}'$ stabilizes $\tilde{L}$ since, given
$\psi\in\Gamma(\tilde{L})$ and using $*A=-AS$, we have
\[
A\wedge\del\psi=A\wedge\hat{\nabla}\psi=-d^{\hat{\nabla}} (A\psi)=0\,.
\]
Therefore, $A$ can be regarded as a holomorphic section of the complex line 
bundle $K\Hom_{+}(\bar{V}/\tilde{L},L)$. If none of the sections $\delta$, $\tilde{\delta}$ 
and $A$ are trivial, we obtain
 \[
\deg K + \deg (V,S) - 2 \deg L\geq 0\,,
\]
\[
\deg K +\deg(V,\tilde{S})-2\deg \tilde{L}\geq 0\,,
\]
\[
\deg K + \deg L - \deg (V,\tilde{S}) + \deg \tilde L\geq 0\,,
\]
and thus
\[
 4\deg K + v -\tilde v\geq 0\,,
\]
with $v$ and $\tilde v$ the degree of the normal bundles of $f$
and $\tilde f$.  Telescoping this last relation over the sequence $f_0,\dots,f_{i-1}$ yields
\[
 v- v_i +4 i \deg K\geq 0\,.
\]
But $4\pi v_i = W(f_i)-W(f_{i-1})$ is a difference of Willmore
energies \eqref{eqn:Will-diff} so that by summing over $i=1,\dots,n$
we get the following relation between the sequence length and the
Willmore energy of $f$:
\[
n\, v +\tfrac{1}{4\pi}W(f)+2n(n+1)\deg K\geq 0\,.
\]
As an immediate consequence we obtain a criterion when the Willmore
sequence has finite length.
\begin{theorem}\label{thm:finite_seq}
  Let $f\colon M\to S^4$ be a branched Willmore surface such that
  either $M$ has genus zero or $M$ has genus one and $f$ has
  non-trivial normal bundle. Then the Willmore sequence of $f$
  terminates.
\end{theorem}
\begin{proof}
  If $M$ has genus zero then $\deg K=-2$ which gives a contradiction
  to the finiteness of $W(f)$. If $M$ has genus one then $\deg K=0$
  and we again contradict the finiteness of $W(f)$, since we always
  may choose our orientation so that the normal bundle degree is
  negative.
\end{proof}

\section{Finite  Willmore sequences}

We have seen that in certain circumstances the Willmore sequence of a
compact branched Willmore surface $f\colon M\to S^4$ terminates. As it
turns out, a finite Willmore sequence can a most have length two and
$f$ is either the twistor projection of a holomorphic curve into
$\CP^3$ or comes from a minimal surface in $\R^4$ by inversion at
infinity.

From the previous section we know that the Willmore sequence can be
continued past a branched Willmore surface $f\colon M\to S^4$ as long
as $A$ and the derivative $\tilde{\delta}$ of $\tilde{f}$ are not
zero.  We first discuss what happens when $\tilde{\delta}=0$ and
hence $\tilde{f}$ is a point, say $\infty\in S^4$. Since the mean
curvature sphere congruence $S$ of $f$ stabilizes $\tilde{L}$ the
spheres $S(p)$ for $p\in M$ all contain the point $\infty$. Therefore,
by inverting $f$ at this point we obtain a surface in $\R^4$ whose
mean curvature spheres become its tangent planes and hence this
surface is minimal.  Since there are no compact minimal surfaces in
$\R^4$ the point $\infty=f(q)$ for some $q\in M$, i.e., $\infty$ must
be lying on the surface $f$.  Moreover, if $f$ is immersed at $q\in M$
then the corresponding minimal surface in $\R^4$ has planar ends.

The fiber $f^{-1}(\infty)\subset M$ is finite because $f$ is a
branched conformal immersion. Thus, away from those points the bundle
$V$ decomposes as a direct sum $V=L\oplus \underline{ L_q}$ of $L$ and
the trivial bundle $\underline{L_q}=M\times L_q$. But $S$ stabilizes
$\underline{L_q}$ and hence also
$Q\underline{L_q}\subset\underline{L_q}$ which, together with $QL=0$,
implies that $\im Q\subset \underline{L_q}$. Therefore, also $\hat{f}$
is the point $\infty$ and the original Willmore sequence consisted of
only one element, $f$, and terminated on either side to the same
constant map. It is easy to see that the converse also holds: a
minimal surface in $\R^4\subset S^4$ is a Willmore surface whose
forward and backward B\"acklund transforms are the point at infinity
$S^4=\R^4\cup \{\infty\}$.

We now come to the second possibility, namely $A=0$, for the Willmore
sequence to terminate. In this case it follows from
\cite[Thm. 4]{coimbra} that $f:M\to S^4$ is the composition of a
holomorphic curve $h:M\to \CP^3$ followed by the twistor projection
$\CP^3\to S^4$.  If in addition $Q=0$, then $dS=0$ by \eqref{eqn:dS},
and $S$ is constant implying that $f(M)\subset S^4$ is itself a round
$2$-sphere. Excluding this possibility we show that the branched
Willmore surface $\hat{f}\colon M\to S^4$ has $\hat{Q}=0$ and is
therefore the twistor projection of a holomorphic curve into $\CP^3$,
after changing the orientation on $S^4$.  To see that $\hat{Q}=0$ we
show that $\hat{S}=-S$ which implies $\hat{Q}=-A=0$. For this it
suffices to observe that
\begin{equation}\label{eqn:delta-hat}
*\hat{\delta}=-\hat{\delta}S=-S\hat{\delta}\,, 
\end{equation}
meaning that $-S$ is a touching sphere congruence along $\hat{f}$. Due
to $A=0$ this touching sphere congruence is in fact the mean curvature
sphere congruence of $\hat{f}$. The first equality in
\eqref{eqn:delta-hat} follows from Theorem~\ref{thm:main} whereas the
second equality follows from the definition of $\hat{\delta}$, the
fact that $A=0$, $\im Q\subset \hat{L}$ and the first equality in
\eqref{eqn:delta-hat}:
\[
\delta=\hat{\pi }d_{|_{\hat{L}}}=\hat{\pi} (\del+\delbar+A+Q)_{|_{\hat{L}}}=\hat{\pi}\delbar_{|_{\hat{L}}}\,. 
\]
Thus, if $f$ does not come from a minimal surface the Willmore
sequence has exactly two elements both of which are twistor
projections of holomorphic curves in $\CP^3$.

Strictly speaking, we have only discussed what happens if the sequence
terminates in the forward direction. To deal with the reversed
direction, we note that we have a dual branched Willmore surface
$f^{\perp}$ obtained from the point-point duality of $\H\P^1$. Its
mean curvature sphere congruence is given by $S^{\perp}=S^*$ and then
$\ker A^{\perp}=(\im Q)^{\perp}$ which reverses the direction of the
Willmore sequence.

This discussion, together with Theorem~\ref{thm:finite_seq}, gives a
unified explanation of results by various authors \cite{Bryant},
\cite{Ejiri}, \cite{montiel}, \cite{coimbra},\cite{willmore_tori} on
the characterization of Willmore surfaces:
\begin{theorem}
  Let $f\colon M\to S^4$ be a branched Willmore surface and assume
  that $M$ has genus zero or that $M$ has genus one and the normal
  bundle of $f$ is non-trivial. Then the Willmore sequence of $f$ has
  at most two elements and $f$ is obtained from a minimal surface in
  $\R^4$ via inversion at infinity, or $f$ is the twistor projection
  of a holomorphic curve in $\CP^3$.
\end{theorem}
 \section{Smooth mean curvature sphere congruences}
 In this last section we provide a proof that the mean curvature
 spheres $\tilde{S}$ of a B\"acklund transform $\tilde{f}\colon M\to
 S^4$ of a branched Willmore surface $f\colon M\to S^4$ extend
 smoothly across the branch locus of $\tilde{f}$, which is the main
 content of Theorem~\ref{thm:main}.  Rather then computing the mean
 curvature spheres for the B\"acklund transform directly, we show that
 an intermediate transform, the $1$-step B\"acklund transform, has
 smooth mean curvature spheres.  We then apply a result of \cite[Lemma
 10]{coimbra} which expresses our B\"acklund transform $\tilde{f}$ as
 a composition of two $1$-step transforms.

First note that $\tilde{f}$ is a conformal map and thus has isolated
branch points: for $\psi\in \Gamma(\tilde{L})$ we have
\[
0=d(*A\psi)=-*A\wedge d\psi = -*A\wedge\tilde{\pi} d\psi=-*A\wedge \tilde{\delta}\psi
\]
where we used $\tilde{L}\subset \ker A$. Since $*A=A(-S)$ and we
assume that $A$ is non-trivial this gives $*\delta=-S\delta$ which
implies conformality of $\tilde{f}$.  Therefore, away from its
isolated branch points $\tilde{f}$ has a mean curvature sphere
congruence $\tilde{S}$ which is given by $-S$ on $V/\tilde{L}$.  It is
shown in \cite[Thm. 7]{coimbra} that $\tilde{Q}=A$ which implies that
$\tilde{S}$ is harmonic and $\tilde{f}$ an immersed Willmore surface
(away from its branch locus). Thus it suffices to show that
$\tilde{S}$ can be extended smoothly across the branch points of
$\tilde{f}$ to obtain a harmonic mean curvature sphere congruence for
$\tilde{f}$ which then is a branched Willmore surface. This is a
purely local problem and so we may assume that $M$ is a simply
connected Riemann surface.

To define a $1$-step B\"acklund transform of the branched Willmore
surface $f$, we choose a point at infinity $\infty=e\H$ in $S^4$ which
is not contained in the image of $f$. After applying a M\"obius
transformation, we may assume $e=\begin{pmatrix}1\\0\end{pmatrix}$ so
that $f=\begin{pmatrix}g\\1\end{pmatrix}\H$.  Let $\alpha\in (\H^2)^*$
be the projection onto the second component so that $<\alpha,e>=0$.
Then, due to the harmonicity \eqref{eqn:harm} of the mean curvature
sphere congruence $S$ of $f$, the $\H$-valued $1$-form
\begin{equation}\label{eqn:BL}
dg^{\sharp}=<\alpha,*Ae>
\end{equation}
is closed and thus exact since $M$ is simply connected.  The map
$g^{\sharp}\colon M\to \H$, respectively the map $f^{\sharp}\colon
M\to S^4$ given by
$f^{\sharp}=\begin{pmatrix}g^{\sharp}\\1\end{pmatrix}\H$, is called a
{\em $1$-step B\"acklund transform} of the branched Willmore surface $f$.
\begin{theorem}\label{lem:1-stepwillmore}
  Let $f\colon M\to S^4$ be a branched Willmore surface and let
  $f^{\sharp}\colon M\to S^4$ be a $1$-step B\"acklund transform. Then
  $f^{\sharp}$ is again a branched Willmore surface, that is to say,
  the mean curvature sphere congruence $S^{\sharp}$ extends smoothly
  across the branch locus of $f^{\sharp}$.
\end{theorem}
\begin{proof}
  Since $\infty\in S^4$ is not on the surface $f$ the trivial bundle
  $V=M\times \H^2$ splits as
\[
V=L\oplus \underline{e\H}
\]
and $V/L=\underline{e\H}$.  We choose the adapted framing $(\psi, e)$
of $V=L\oplus\underline{e\H}$ with
$\psi=\begin{pmatrix}g\\1\end{pmatrix}$ a smooth section of $L$. In
this frame the mean curvature sphere congruence $S$ of $f$ is given by
\cite[Sec. 7.2]{coimbra}
\[
S\psi=-\psi R\quad\text{and}\quad Se=eN -\psi H 
\]
with $R,N,H\colon M \to\H$ satisfying $R^2=N^2=-1$ and $RH=NR$. The
conformality and touching conditions $*\delta=S\delta=\delta S$
translate to
\[
*dg=N dg= -dg R
\]
so that $N$ and $R$ are the left and right normals of $g:M\to \H$.
That $S$ is the mean curvature sphere congruence along $f$ is
expressed by \cite[Sec. 7.2]{coimbra}
\begin{equation}\label{eqn:mcs}
2dg H= dN-N*dN\,.
\end{equation}
{F}rom \cite[Thm. 6]{coimbra} we know that the $1$-step B\"acklund transform  $f^{\sharp}$ is a branched conformal immersion which is Willmore away from its branch points.  Therefore 
it suffices to show that the mean curvature sphere congruence $S^\sharp$ of  
$f^{\sharp}$ extends across the branch points which is a purely local argument.  To calculate the left
and right normals for the $1$-step transform $g^{\sharp}$, we need to
split the trivial bundle $V$ with respect to the B\"acklund transform
$\tilde{f}$ into $V=\tilde{L}\oplus \underline{e\H}$.  This is
possible since the point $\infty\in S^4$ can be chosen
 in the complement of the images of $f$ and $\tilde{f}$, after possibly
restricting to a sufficiently small neighborhood of a branch point of
$f^\sharp$.
Let $\beta\in\Gamma(V^{*})$ be the unique smooth section satisfying
$<\beta,e>=1$ and $\beta_{|\tilde{L}}=0$.  We now express
$Se=\tilde{\phi}+e<\beta,Se>$ in the splitting $V=\tilde{L}\oplus
\underline{e\H}$ and recall that $\tilde{L}\subset \ker A$ and
$*A=-AS$. Then, from the definition \eqref{eqn:BL} of $g^{\sharp}$, we
see that its right normal is given by
\[
\psi *dg^{\sharp} = S*Ae=-*ASe=-*Ae<\beta,Se>=-\psi dg^{\sharp} <\beta,Se>\,.
\]
On the other hand, from $*A=SA$ we obtain 
\[
\psi *dg^{\sharp}=-Ae=S\psi dg^{\sharp}=-\psi R dg^{\sharp}\,.
\]
Therefore, the left and right normals of $g^{\sharp}$ are given by
\[
N^{\sharp}=-R\quad\text{and}\quad R^{\sharp}=<\beta,Se>
\]
which are both smooth on $M$.  Hence, if the mean curvature sphere
congruence $S^{\sharp}$ exists it is expressible in the adapted frame
$(\psi^{\sharp},e)$ of $V=L^{\sharp}\oplus \underline{e\H}$ by
\[
S^{\sharp}\psi^{\sharp}=-\psi^{\sharp}R^{\sharp}\quad\text{and}\quad S^{\sharp}e=eN^{\sharp}- \psi^{\sharp}H^{\sharp}\,,
\]
with $\psi^{\sharp}=\begin{pmatrix}g^{\sharp}\\1\end{pmatrix}$ a
smooth section of $L^{\sharp}$.  In particular, $S^{\sharp}$ has to
satisfy
\begin{equation}\label{eqn:mean-sharp}
2dg^{\sharp}H^{\sharp}=dN^{\sharp}-N^{\sharp}*dN^{\sharp}\,
\end{equation}
by \eqref{eqn:mcs}, which can be used to calculate $H^{\sharp}$: since
$N^{\sharp}=-R$ we have
\[
dN^{\sharp}-N^{\sharp}*dN^{\sharp}=-dR-R*dR=-4<\alpha,*A\psi>\,,
\]
with the last equality following from \cite[Prop. 12]{coimbra}.
Decomposing $\psi=\tilde{\phi}+e<\beta,\psi>$ in the splitting
$V=\tilde{L}\oplus\underline{e\H}$, and recalling once more
$\tilde{L}\subset \ker A$ and the definition \eqref{eqn:BL} of
$dg^{\sharp}$, we obtain from \eqref{eqn:mean-sharp}
\[
dg^{\sharp}(H^{\sharp}+2<\beta,\psi>)=0
\]
and hence $H^{\sharp}=-2<\beta,\psi>$ is smooth on $M$.  For
$S^{\sharp}$ to actually be a sphere congruence we need
$(S^{\sharp})^2=-1$. But $(R^{\sharp})^2=(N^{\sharp})^2=-1$ so it
remains to verify $R^{\sharp}H^{\sharp}=H^{\sharp}N^{\sharp}$ which
immediately follows from the explicit expression of $H^{\sharp}$.
\end{proof}

As already mentioned it is proven in \cite[Lemma 10]{coimbra} that the
B\"acklund transform $\tilde{f}$ of $f$ is given by two successive
1-step transforms. Thus, applying the previous theorem twice we see
that the mean curvature sphere congruence $\tilde{S}$ of $\tilde{f}$
is smooth.  To obtain the result for the B\"acklund transform
$\hat{f}$, we apply the above argument to the dual surface $f^{\perp}$ thereby concluding
the proof of Theorem~\ref{thm:main}.
\begin{rem}
  There is a more conceptual way to prove
  Theorem~\ref{lem:1-stepwillmore}: using the theory of envelopes of
  Frenet curves \cite{osculates} in $\HP^n$ one can show that a 1-step
  B\"acklund transform is, up to a suitable projection, given by an
  envelope \cite{onestep}. A Frenet curve is the natural
  generalization to $\HP^n$ of a conformal map into $S^4$ allowing a
  smooth mean curvature sphere congruence and the enveloping
  construction preserves the Frenet property.  Using this generalized
  setup an analogue of Theorem~\ref{thm:main} for Willmore surfaces in $\H\P^n$ can be found in
  \cite{habil}.
\end{rem}

\bibliographystyle{plain}
\bibliography{doc}  
\end{document}